\documentclass{amsart}
\usepackage{amssymb,latexsym,graphics,amsfonts}

\swapnumbers

\theoremstyle{definition}

\theoremstyle{definition}

\numberwithin{equation}{section}
 \numberwithin{equation}{subsection}

\begin{document}

\title{module super-amenability for semigroup algebras}

\author[A. Bodaghi]{Abasalt Bodaghi}
\address{Department of Mathematics, Islamic Azad University, Garmsar Branch, Garmsar,
 Iran}
\email{abasalt\_bodaghi@yahoo.com, abasalt.bodaghi@gmail.com}

\author[M. Amini]{Massoud Amini}
\address{Department of Mathematics,
 Tarbiat Modares University, Tehran 14115-175, Iran}

\email{mamini@modares.ac.ir}

\subjclass[2000]{46H25}

\keywords{Banach modules, module derivation, module
super-amenability, inverse semigroups}

\dedicatory{}

\smallskip

\begin{abstract}
Let $S$ be an inverse semigroup with the set of idempotents $E$.
In this paper we define the module super-amenability of a Banach
algebra which is a Banach module over another Banach algebra with
compatible actions, and show that when $E$ is upward directed and
acts on $S$ trivially from left and by multiplication from right,
the semigroup algebra $ \ell ^{1}(S)$ is $\ell^{1}(E)$-module
super-amenable if and only if an appropriate group homomorphic
image of $S$ is finite.
\end{abstract}

\maketitle

\section{introduction}
The second author in \cite{am1} introduced the concept of module
amenability and showed that for an inverse semigroup $S$, the
semigroup algebra $\ell^1(S)$ is module amenable as a Banach
module on $\ell^1(E)$, where $E$ is the set of idempotents of
$S$, if and only if $S$ is amenable (see also [2]). This is the
semigroup analog of Johnson's theorem for locally compact groups
\cite{joh}. In this paper, we find a similar result for module
super-amenability of the semigroup algebra of an inverse
semigroup which is the semigroup analog of Selivanov's theorem
for locally compact groups \cite{sel}.

Recall that a Banach algebra ${\mathcal A}$ is called {\it
super-amenable} ({\it contractible}) if $H^1({\mathcal A},\mathcal
X)=\{0\}$ for every Banach ${\mathcal A}$-bimodule $\mathcal X$,
where the left hand side is the {\it first cohomology group} of
${\mathcal A}$ with coefficient in $\mathcal X$ (see \cite{dal,
run}). A Banach space $E$ has the {\it approximation property} if
there is a net $(T_j)_j$ in $\mathcal F(E)$, the space of the
bounded finite rank operators on $E$ such that $T_j\longrightarrow
id_E$ uniformly on compact subsets on $E$. It is shown in
\cite[Theorem 4.1.5]{run} if $\mathcal A$ is a super-amenable Banach
algebra and has the approximation property, then $\mathcal A$ is
finite dimensional (see also Propositions 5.1 and 5.2 of
\cite{tay}). In particular, since $\ell^{1}(S)$ has the
approximation property \cite{dal}, so it is not super-amenable, when
$S$ is infinite. For groups, Selivanov showed in \cite{sel} that for
any locally compact group $G$, $L^1(G)$ is super-amenable if and
only if $G$ is finite (see also \cite[Exercise 4.1.3]{run}). In this
paper, we develop the concept of {\it module super-amenability} for
a class of Banach algebras, and prove that for an inverse semigroup
$S$ with an upward directed subsemigroup of idempotents $E$ which
acts on $S$ trivially from left and by multiplication from right,
the semigroup algebra $ \ell ^{1}(S)$ is $\ell^{1}(E)$-module
super-amenable if and only if an appropriate group homomorphic image
of $S$ is finite.

The paper is organized as follows: Section 2 is devoted to the
concept of module super-amenability. The main result of this
section asserts that module super-amenability is equivalent to
the existence of a module diagonal. In section 3 we show that for
an inverse semigroup $S$ with an upward directed set of
idempotents $E$, $\ell^1(S)$ is $\ell^1(E)$-module super-amenable
(with respect to a the above action) if and only if the group
homomorphic image $S/\approx$ is finite, where $s\approx t$
whenever $\delta_s-\delta_t$ belongs to the closed linear span of
the set
$$\{\delta_{set}-\delta_{st}: s,t\in S, e\in E\}.$$
Examples of semigroups with an upward directed set of idempotents
include all unital inverse semigroups, the bicyclic semigroup,
and the semigroup of natural numbers with {\it max} operation.
The free inverse semigroup on two generators is an example of an
inverse semigroup whose set of idempotents is not upward directed.

\section{Module Super Amenability  }

Throughout this paper, ${\mathcal A}$ and ${\mathfrak A}$ are
Banach algebras such that ${\mathcal A}$ is a Banach ${\mathfrak
A}$-bimodule with compatible actions, that is

$$
\alpha\cdot(ab)=(\alpha\cdot a)b,
\,\,(ab)\cdot\alpha=a(b\cdot\alpha) \hspace{0.3cm}(a,b \in
{\mathcal A},\alpha\in {\mathfrak A}).
$$
Let ${\mathcal X}$ be a Banach ${\mathcal A}$-bimodule and a
 Banach ${\mathfrak A}$-bimodule with compatible actions, that is
$$
\alpha\cdot(a\cdot x)=(\alpha\cdot a)\cdot x,
\,\,a\cdot(\alpha\cdot x)=(a\cdot\alpha)\cdot x, \,\,(\alpha\cdot
x)\cdot a=\alpha\cdot(x\cdot a) \hspace{0.3cm}(a \in{\mathcal
A},\alpha\in {\mathfrak A},x\in{\mathcal X} )
$$
and the same for the right or two-sided actions. Then we say that
${\mathcal X}$ is a Banach ${\mathcal A}$-${\mathfrak A}$-module.
If moreover
$$\alpha\cdot x=x\cdot\alpha \hspace{0.3cm}( \alpha\in {\mathfrak
A},x\in{\mathcal X} )$$ then $\mathcal X $ is called a {\it
commutative} ${\mathcal A}$-${\mathfrak A}$-module. If $\mathcal X
$ is a (commutative) Banach ${\mathcal A}$-${\mathfrak
A}$-module, then so is $\mathcal X^*$, where the actions of
$\mathcal A$ and ${\mathfrak A}$ on $\mathcal X^*$ are defined by
$$\langle\alpha\cdot f,x\rangle{}=\langle{}f,x\cdot\alpha\rangle{},\,\,\langle{}
a\cdot f,x\rangle{}=\langle{}f,x\cdot a\rangle{}\hspace{0.3cm} (a
\in{\mathcal A},\alpha\in {\mathfrak A},x\in{\mathcal X},f \in
\mathcal X^* )$$
 and the same for the right actions. Let ${\mathcal Y}$ be another ${\mathcal A}$-${\mathfrak
A}$-module, then a ${\mathcal A}$-${\mathfrak A}$-module morphism
from ${\mathcal X}$ to ${\mathcal Y}$ is a norm-continuous map
$\varphi :{\mathcal X}\longrightarrow {\mathcal Y}$ with $
\varphi(x \pm y)= \varphi (x)\pm \varphi (y)$ and
 $$ \varphi (\alpha \cdot x)= \alpha \cdot\varphi(x), \,\,\varphi(x\cdot \alpha)= \varphi (x)\cdot\alpha
,\,\,\varphi (a\cdot x)=a\cdot \varphi(x), \varphi(x\cdot a)=
\varphi(x)\cdot a,$$ for $x,y \in {\mathcal X}, a \in {\mathcal
A},$ and $\alpha \in {\mathfrak A}$.

 Note that when
${\mathcal A}$ acts on itself by algebra multiplication, it is
not in general a Banach ${\mathcal A}$-${\mathfrak A}$-module, as
we have not assumed the compatibility condition
$$a\cdot(\alpha\cdot b)=(a\cdot\alpha)\cdot b\quad (\alpha\in {\mathfrak A}, a,b \in{\mathcal
A}).$$ If $\mathcal A$ is a commutative $\mathfrak A$-module and
acts on itself by multiplication from both sides, then it is also
a Banach ${\mathcal A}$-${\mathfrak A}$-module.

If ${\mathcal A}$ is a Banach $\mathfrak A$-module with
compatible actions, then so are the dual space $\mathcal A^*$ and
the second dual space $\mathcal A^{**}$. If moreover $\mathcal A$
is a commutative $\mathfrak A$-module, then $\mathcal A^*$ and the
$\mathcal A^{**}$ are commutative ${\mathcal A}$-${\mathfrak
A}$-modules.

Consider the projective tensor product $\mathcal A \widehat
\bigotimes \mathcal A$. It is well known that $\mathcal A
\widehat \bigotimes \mathcal A$ is a Banach algebra with respect
to the canonical multiplication map defined by
$$(a\otimes b)(c\otimes d)=(ac\otimes bd)$$
and extended by bi-linearity and continuity [4]. Then $\mathcal A
\widehat \bigotimes \mathcal A$ is a Banach ${\mathcal
A}$-${\mathfrak A}$-module with canonical actions. Let $I$ be the
closed ideal of the projective tensor product $\mathcal A
\widehat \bigotimes \mathcal A$ generated by elements of the form
$\alpha \cdot a \otimes b-a \otimes b\cdot\alpha$ for $ \alpha\in
{\mathfrak A},a,b\in{\mathcal A}$. Consider the map $\omega \in
{\mathcal L}(\mathcal A \widehat \bigotimes \mathcal A, \mathcal
A )$ defined by $\omega (a \otimes b)=ab$ and extended by
linearity and continuity. Let $J$ be the closed ideal of
${\mathcal A}$ generated by  $\omega(I)$. Then the module
projective tensor product ${\mathcal A}\widehat \bigotimes
_{\mathfrak A} {\mathcal A}\cong(\mathcal A \widehat \bigotimes
\mathcal A)/{I}$ and the quotient Banach algebra $\mathcal A/J$
are Banach ${\mathfrak A}$-modules with compatible actions. Also
the map $\widetilde{\omega} \in {\mathcal L}({\mathcal A}\widehat
\bigotimes _{\mathfrak A} {\mathcal A},\mathcal A/J)$ defined by
$ \widetilde{\omega} (a \otimes b +I)=ab+J$ extends to an
${\mathfrak A}$-module morphism.

Let ${\mathcal A}$ and ${\mathfrak A}$ be as above and ${\mathcal
X}$ be a Banach ${\mathcal A}$-${\mathfrak A}$-module. A bounded
map $D: \mathcal A \longrightarrow \mathcal X $ is called a {\it
module derivation} if
$$D(a\pm b)=D(a)\pm D(b),\hspace{0.2cm}D(ab)=D(a)\cdot b+a\cdot D(b)\hspace{0.3cm}
(a,b \in \mathcal A),$$and
$$D(\alpha\cdot a)=\alpha\cdot D(a),\hspace{0.3cm}D(a\cdot\alpha)=D(a)\cdot\alpha
\hspace{0.3cm}(a \in{\mathcal A},\alpha\in {\mathfrak A}).$$ Note
that $D: \mathcal A \longrightarrow \mathcal X $ is bounded if
there exist $M>0$ such that $\| D(a) \| \leq M \| a \| $, for
each $a \in{\mathcal A}$. Although $D$ is not necessarily linear,
but still its boundedness implies its norm continuity (since $D$
preserves subtraction). When $\mathcal X $ is commutative, each
$x \in \mathcal X $ defines a module derivation
$$D_x(a)=a\cdot x-x\cdot a \hspace{0.3cm} (a \in{\mathcal A}).$$
These are called {\it inner} module derivations.

\vspace{.3cm}\paragraph{\large\bf Definition 2.1.}  {\it The
Banach algebra ${\mathcal A}$ is called module super-amenable (as
an ${\mathfrak A}$-module) if for any commutative Banach
${\mathcal A}$-${\mathfrak A}$-module $\mathcal X$, each module
derivation $D: \mathcal A \longrightarrow \mathcal X$ is inner.}

\vspace{.3cm} We use the notations $Z_{\mathfrak A}(\mathcal A,
\mathcal X)$ and $B_{\mathfrak A}(\mathcal A, \mathcal X)$ for the
set of all module derivations and inner module derivations from
$\mathcal A$ to $\mathcal X$, respectively. The quotient group
(called the first relative -to ${\mathfrak A}$- cohomology group
of ${\mathcal A}$ with coefficients in $\mathcal X$) is denoted
by $H^1_{\mathfrak A}(\mathcal A, \mathcal X)$. Hence ${\mathcal
A}$ is module super-amenable if and only if $H^1_{\mathfrak
A}(\mathcal A, \mathcal X)=\{0\}$, for each commutative Banach
${\mathcal A}$-${\mathfrak A}$-module $\mathcal X$.

\vspace{.3cm}\paragraph{\large\bf Proposition 2.2.} {\it If
${\mathfrak A}$ has an identity for $\mathcal A$, then
super-amenability of $\mathcal A$ implies its module
super-amenability.}

\vspace{.3cm}\paragraph{\large\bf Proof.} Let ${\mathfrak e}\in
\mathfrak A$ be a identity for $\mathcal A$, that is $\mathfrak
e. a=a.\mathfrak e=a$, for each $a\in\mathcal A$, and $\mathcal
X$ be a commutative ${\mathcal A}$-${\mathfrak A}$-module. Assume
that $D:\mathcal A \longrightarrow \mathcal X$ is a module
derivation, then obviously $D(a\cdot\lambda \mathfrak
e)=D(\lambda a)$, for each $a \in{\mathcal A}$ and $\lambda \in
\mathbb{C}$. On the other hand,
$$D(a\cdot\lambda \mathfrak e) = D(a)\cdot\lambda \mathfrak e
= \lambda D(a)\cdot\mathfrak e= \lambda D(a\cdot\mathfrak
e)=\lambda D(a).$$ Thus $D$ is $\mathbb{C}$-linear, and so
inner.$\hfill\square $

\vspace{.3cm} As we will see later in section 3, there are module
super-amenable Banach algebras that are not super-amenable, so
the converse of the above Proposition is false. It is known that
every super-amenable Banach algebra has an identity. Also a
Banach algebra is super-amenable if and only if it has a diagonal
\cite{run}. Recall that a {\it diagonal} for $\mathcal A$ is an
element $M\in \mathcal A\widehat\bigotimes \mathcal A$ satisfying
$$a\cdot\omega ( M)=a,\hspace{0.2cm}a\cdot M=M\cdot a\hspace{0.2cm}(a\in
\mathcal A).$$ We start this section by showing that similar
results hold $\mathcal A$ is a commutative $\mathcal
A$-$\mathfrak A$-module.

\vspace{.3cm}\paragraph{\large\bf Proposition 2.3.} {\it Let
$\mathcal A$ be a commutative Banach $\mathcal A$-${\mathfrak
A}$-module. If $\mathcal A$ is module super-amenable, then it is
unital.} \vspace{.3cm}\paragraph{\large\bf Proof.} Let's consider
$\mathcal X=\mathcal A$ as an $\mathcal A$-bimodule, with actions
$$ a\cdot b:=ab,\hspace{0.4cm} b\cdot a:=0 \hspace{0.3cm} (a\in \mathcal
A, b\in\mathcal X).$$

Let $D:\mathcal A\longrightarrow \mathcal X$ be the identity map,
it is clear that $D \in Z_{\mathfrak A}(\mathcal A, \mathcal
X)=B_{\mathfrak A}(\mathcal A,\mathcal X)$. This means that there
is $a_0 \in \mathcal A$ such that $aa_0=a$, for all $a \in
\mathcal A.$ Therefore $a_0$ is a right identity for $\mathcal
A$. Similarly $\mathcal A$ has a left identity. The left and
right identities now have to coincide.$\hfill\square $

\vspace{.3cm}\paragraph{\large\bf Proposition 2.4.} {\it Let
$\mathcal A$ and $\mathcal B$ be Banach algebras and Banach
${\mathfrak A}$-modules with compatible actions. If $\mathcal A$
is ${\mathfrak A}$-module super-amenable and $\varphi:\mathcal A
\longrightarrow \mathcal B$ is a continuous Banach algebra
homomorphism with dense range, then $\mathcal B$ is also
${\mathfrak A}$-module super-amenable.}

\vspace{.3cm}\paragraph{\large\bf Proof.} Let $X$ be a
commutative ${\mathcal B}$-${\mathfrak A}$-module, then it is a
commutative ${\mathcal A}$-${\mathfrak A}$-module with the
following actions
$$a\cdot x:=\varphi(a)\cdot x,\hspace{0.4cm}x\cdot a:=x\cdot\varphi(a) \hspace{0.3cm}(a \in \mathcal
A, x\in X).$$

If $D:\mathcal B \longrightarrow X$ is a module derivation, then
$D\circ \varphi:\mathcal A \longrightarrow X$ is a module
derivation, which is inner. By density of the range of $\varphi$
and continuity of $D$, $D$ is inner.$\hfill\square $

\vspace{.3cm}\paragraph{\large\bf Definition 2.5.} {\it An element
$M \in {\mathcal A}\widehat \bigotimes _{\mathfrak A} {\mathcal
A}$ is called a module diagonal if $\widetilde {\omega}(M)$ is an
identity of $\mathcal A/J$ and $a\cdot M=M\cdot a$, for all $a\in
\mathcal A$.}

\vspace{.3cm}\paragraph{\large\bf Theorem 2.6.} {\it Let
$\mathcal A$ be a Banach $\mathcal A$-${\mathfrak A}$-module.
Then $\mathcal A$ is module super-amenable if and only if
$\mathcal A$ has a module diagonal.}

\vspace{.3cm}\paragraph{\large\bf Proof.} Assume that $\mathcal A$
is module super-amenable, then by Proposition 2.3 it has an
identity element $e$. Put $T=e\otimes e+I$, we have
$\widetilde{\omega}(a\cdot T-T\cdot a)=J$. Hence
$\widetilde{\omega}$ vanishes on the range of $D_T$, and $D_T$
could be regarded as a module derivation into $K$. Since
${\mathcal A}\widehat \bigotimes _{\mathfrak A} {\mathcal A}$ is
always commutative $\mathfrak A$-module, so is $K$=ker
$\widetilde{\omega}$, hence by module super-amenability of
${\mathcal A}$, there is $N\in K$ such that $D_T=D_N$. Now it is
easy to see that $M=N-T$ is a module diagonal.

Conversely suppose that $M=\sum_{n=1}^{\infty} a_n\otimes b_n +I$
is a module diagonal, where $(a_n),(b_n)$ are bounded sequences
in $\mathcal A$ with $\sum_{n=1}^{\infty} \|a_n\|\|
b_n\|<\infty$. Let $\mathcal X$ be a commutative Banach
${\mathcal A}$-${\mathfrak A}$-module, then clearly $J$ acts
trivially on $\mathcal X$, that is $J\cdot{\mathcal X}={\mathcal
X}\cdot J=\{0\}$. Therefore $\mathcal X$ is a Banach ${\mathcal
A}/{J}$-module with the following actions
$$(a+J)\cdot x:=a\cdot x, \hspace{0.2cm}x\cdot(a+J):=x\cdot a \hspace{0.3cm} (x \in \mathcal X ,a \in
\mathcal A).$$

If $D:\mathcal A \longrightarrow \mathcal X$ is a module
derivation, then the map $\tilde{D}:\mathcal A/J \longrightarrow
\mathcal X$ defined by $\tilde{D}(a+J)= D(a)$, for $a \in
\mathcal A$, is a module derivation. Consider
$x=\sum_{n=1}^{\infty} a_n\cdot\tilde{D}(b_n+J)$, then for each
$\phi \in \mathcal X^*$ we have
\begin{align*}
\langle{}\phi,(a+J)\cdot x \rangle{}
&=\langle{}\phi,(a+J)\cdot\sum_{n=1}^{\infty}a_n\cdot\tilde{D}(b_n+J)\rangle{} \\
&=\langle{}\phi,\sum_{n=1}^{\infty}a_n\cdot\tilde{D}(b_na+J)\rangle{} \\
&=\langle{}\phi,\sum_{n=1}^{\infty}a_n\cdot\tilde{D}(b_n+J)\cdot(a+J)\rangle{}+\langle{}\phi,(\sum_{n=1}^{\infty}
a_nb_n+J)\cdot\tilde{D}(a+J)\rangle{}  \\
&=\langle{}\phi,x\cdot(a+J)\rangle{}+\langle{}\phi,\tilde{\omega}(M)\cdot\tilde{D}(a+J)\rangle{}\\
&=\langle{}\phi,x\cdot(a+J)\rangle{}+\langle{}\phi,\tilde{D}(a+J)\rangle{}.
\end{align*}
Hence $\tilde{D}$ is inner. Therefore $D$ is an inner module
derivation. $\hfill\square $

 We say the Banach algebra ${\mathfrak A}$ acts trivially on $\mathcal A$ from left if for
each $\alpha\in \mathfrak A$ and $a\in \mathcal A$, $\alpha\cdot
a=f(\alpha)a$, where $f$ is a continuous linear functional on
${\mathfrak A}$ (see \cite[Lemma 3.1]{abe}).

 \vspace{.3cm}\paragraph{\large\bf Lemma 2.7.}
{\it Let $\mathcal A$ be a Banach ${\mathfrak A}$-module with
trivial left action. If $\mathcal A$ is module super-amenable and
$\mathcal A/J$ is commutative ${\mathfrak A}$-module, then
$\mathcal A/J$ is super-amenable.}

\vspace{.3cm}\paragraph{\large\bf Proof.} Let $\mathcal X$ be a
commutative $\mathcal A/J$-${\mathfrak A}$-module. Then $\mathcal
X$ is a commutative ${\mathcal A}$-${\mathfrak A}$-module with
the same actions over ${\mathfrak A}$ and module actions over
${\mathcal A}$ defined by
$$a\cdot x:=(a+J)\cdot x,\hspace{0.2cm}x\cdot a:=x\cdot(a+J) \hspace{0.3cm} (x \in \mathcal X ,a \in
\mathcal A).$$ Suppose that $D: \mathcal A/J\longrightarrow
\mathcal X $ is a module derivation, then $\tilde{D}: \mathcal A
\longrightarrow \mathcal X$ defined by $\tilde{D}(a)=D(a+J)$, for
$a \in \mathcal A$, is a module derivation. Since ${\mathcal A}$
is module super-amenable, $\tilde{D}$, and so $D$ are both inner.
Hence $\mathcal A/J$ is module super-amenable, and since
$\mathcal A/J$ is a commutative $\mathfrak A$-module, by
Proposition 2.3, it has an identity. The rest of the proof goes
exactly like that of \cite[Theorem 3.2]{abe}, leading to
super-amenability of $\mathcal A/J$.$\hfill\square $

\vspace{.3cm}\paragraph{\large\bf Lemma 2.8.} {\it Let $\mathcal
A$ be an ${\mathfrak A}$-module with trivial left action. If
${\mathfrak A}$ has a bounded approximate identity for $\mathcal
A$ and $\mathcal A/J$ is commutative ${\mathfrak A}$-module, then
super-amenability of $\mathcal A/J$ implies module
super-amenability of $\mathcal A$.}

\vspace{.3cm}\paragraph{\large\bf Proof.} Let ${\mathcal X}$ be a
commutative Banach ${\mathcal A}$-${\mathfrak A}$-module and $D:
\mathcal A \longrightarrow {\mathcal X} $ be a module derivation.
Since $J\cdot{\mathcal X}={\mathcal X}\cdot J=0$, ${\mathcal X}$
is a Banach $\mathcal A/J$-module with module actions
$$(a+J)\cdot x:=a\cdot x, \hspace{0.2cm}x\cdot(a+J):=x\cdot a \hspace{0.3cm} (x \in \mathcal X ,a \in
\mathcal A).$$

Consider $\tilde{D}: \mathcal A/J\longrightarrow \mathcal X $,
defined by $\tilde{D}(a+J)= D(a)$, for $a \in
 \mathcal A$.  For each $a,b \in
 \mathcal A$ and $\alpha\in \mathfrak A$ we have
$$
D(\alpha\cdot ab-ab\cdot\alpha)=\alpha\cdot
D(ab)-D(ab)\cdot\alpha=0.$$

 By the above observation, $\tilde{D}$ is
also well-defined. Since $\mathcal A/J$ is a commutative
${\mathfrak A}$-module, its super-amenability implies that it has
identity $e+J$. Without loss of generality, we may assume that
$\mathcal X$ is also a unital $\mathcal A/J$-bimodule. Now
${\mathfrak A}$ acts on $\mathcal A$ trivially from left,
 hence $f(\alpha)a-a\cdot\alpha\in J$, for each $\alpha \in {\mathfrak A}$, where
 $f$ is a continuous linear functional on ${\mathfrak A}$ \cite[Lemma 3.1]{abe}.
 Suppose that ${\mathfrak A}$ has a bounded approximate identity $(\gamma_i)$
for $\mathcal A$. Since $f$ is bounded, $\{|f(\gamma_i)|\}$ is a
bounded net in $\mathbb{C}$. Without loss of generality, we may
assume that $f(\gamma_i)\longrightarrow 1$. Thus for each
$\lambda \in \mathbb{C}$ we have
$$e\cdot(\lambda \gamma_i)-f(\gamma_i)e= (\lambda e)\cdot \gamma_i-f(\gamma_i)e
\longrightarrow \lambda e-e$$ in norm. Since $J$ is a closed
ideal of $\mathcal A$, $\lambda e-e \in J$. Next, for $\lambda
\in \mathbb{C}, a\in \mathcal A$, we have
\begin{align*}
\tilde{D}(\lambda(a+J)) &=\tilde{D}((a+J)(\lambda e+J)) \\
&=(a+J)\cdot\tilde{D}(\lambda e+J)+\tilde{D}(a+J)(\lambda e+J)\\
&=(a+J)\cdot \tilde{D}(e+J)+ \lambda \tilde{D}(a+J)\cdot(e+J)\\
&=\lambda \tilde{D}(a+J).
\end{align*}
Thus $\tilde{D}$ is $\mathbb{C}$-linear, and so it is inner.
Therefore $D$ is an inner module derivation.$\hfill\square $

\section{Module Super-Amenability for Semigroup Algebras  }
In this section we find conditions on a (discrete) inverse
semigroup $S$ such that the semigroup algebra $\ell^1(S)$ is
$\ell^1(E)$-module super-amenable, where $E$ is the set of
idempotents of $S$, acting naturally on it. We start this section
with the definition of inverse semigroups.

\vspace{.3cm}\paragraph{\large\bf Definition 3.1.} {\it A
discrete semigroup $S$ is called an inverse semigroup if for each
$ s \in S $ there is a unique element $s^* \in S$ such that
$ss^*s=s$ and $s^*ss^*=s^*$. An element $e \in S$ is called an
idempotent if $e=e^*=e^2$. The set of idempotents of $S$ is
denoted by $E$.}

\vspace{.3cm} There is a natural order on $E$, defined by
$$e\leq d \Longleftrightarrow ed=e \hspace{0.3cm}(e,d \in
E).$$ It is easy to see that $E$ is indeed a commutative
subsemigroup of $S$. In particular $ \ell ^{1}(E)$ could be regard
as a subalgebra of $ \ell ^{1}(S)$, and thereby $ \ell ^{1}(S)$
is a Banach algebra and a Banach $ \ell ^{1}(E)$-module with
compatible canonical actions. However, for technical reasons,
here we let $ \ell ^{1}(E)$ act on $ \ell ^{1}(S)$ by
multiplication from right and trivially from left, that is
$$\delta_e\cdot\delta_s = \delta_s, \,\,\delta_s\cdot\delta_e = \delta_{se} =
\delta_s * \delta_e \hspace{0.3cm}(s \in S,  e \in E).$$

In this case, $J$ is the closed linear span of
$$\{\delta_{set}-\delta_{st} \hspace{0.2cm} s,t \in S,  e \in E
\}.$$ We consider the following equivalence relation on $S$
$$s\approx t \Longleftrightarrow \delta_s-\delta_t \in J \hspace{0.2cm} (s,t \in
S).$$

Recall that $E$ is called {\it upward directed} if for every $e,f
\in E$ there exist $g \in E$ such that $eg=e$ and $fg=f$. This is
precisely the assertion that $S$ satisfies the condition $D_1$ of
Duncan and Namioka \cite{dun}. It is shown in \cite{abe} that if
$E$ is upward directed, then the quotient $S/\approx$ is a
discrete group. Unital inverse semigroups have an upward directed
set of idempotents. Also if $E$ is totally ordered, it is clearly
upward directed. The examples of the latter include the bicyclic
semigroup and the semigroup of natural numbers with {\it max}
operation. On the other hand, the set of idempotents of the free
inverse semigroup on two generators is not upward directed.
Indeed, if the generators are $a$ and $b$, there is no idempotent
which is bigger than both $aa^*$ and $bb^*$.

With the notations of previous section, $ \ell ^{1}(S)/{J}\cong
{\ell ^{1}}(S/\approx)$ is a commutative $ \ell ^{1}(E)$-bimodule
with the following actions
$$\delta_e\cdot(\delta_s+J) = \delta_s+J, \,\,(\delta_s+J)\cdot\delta_e = \delta_{se}+J \hspace{0.3cm}(s \in S,  e \in E).$$

The main theorem of this section is a semigroup analog of the
Selivanov's theorem \cite{sel} for groups, characterizing module
super-amenability of the semigroup algebra of an inverse
semigroup with an upward directed set of idempotents. Indeed we
reduce the result for inverse semigroups to that of discrete
groups, and use Selivanov's theorem.

\vspace{.3cm}\paragraph{\large\bf Theorem 3.2.} {\it Let $S$ be an
inverse semigroup with an upward directed set of idempotents $E$.
 Then $\ell ^{1}(S)$ is module super-amenable, as an $\ell
^{1}(E)$-module with trivial left action and canonical right
action, if and only if $S/\approx$ is finite.}

\vspace{.3cm}\paragraph{\large\bf Proof.} Suppose that $\ell
^{1}(S)$ is module super-amenable, then ${\ell ^{1}(S)}/{J}\cong
{\ell ^{1}}(S/\approx)$ is super-amenable by Lemma 2.7. Since
$S/\approx$ is a (discrete) group, it has to be finite by
Selivanov's theorem \cite{sel}. Conversely, if $S/\approx$ is
finite, then $ {\ell ^{1}(S)}/{J}$ is super-amenable \cite{sel}.
Since $E$ satisfies condition $D_1$ of Duncan and Namioka, so
$\ell ^{1}(E)$ has a bounded approximate identity for $\ell
^{1}(S)$ \cite{abe, dun}. Now the result follows from Lemma 2.8
with $\mathcal A= \ell ^{1}(S)$ and $\mathfrak A=\ell ^{1}(E)$.
$\hfill\square $

\vspace{.3cm} We close this section by some examples of module
super-amenable Banach algebras. Let $\mathfrak G$ be a
commutative unital Banach algebra with unit element $\mathfrak
e$. Consider $\mathbb{A}=M_n(\mathfrak G)$, the Banach algebra of
$n\times n$ matrices with entries from $\mathfrak G$. Then
$\mathbb{A}$ is a unital commutative $\mathfrak G$-bimodule with
the following natural actions
$$\alpha\cdot[\beta_{ij}]=[\alpha \beta_{ij}],\hspace{0.2cm}[\beta_{ij}]\cdot\alpha=[\beta_{ij} \alpha]\hspace{0.2cm}
(\alpha\in \mathfrak G, [\beta_{ij}] \in \mathbb{A}).$$ Consider
the set of matrix units $\{E_{ij} ; i,j=1,...,n \}$, where
$E_{ij}$ is the matrix having $\mathfrak e$ at the $i^{th}$ row
and $j^{th}$ column, and zero elsewhere. The identity matrix $E$,
which is the unit element of $\mathbb{A}$, is the matrix whose
diagonal entries are $\mathfrak e$ and has zero entries
elsewhere. Let $I,J$ be the corresponding closed ideals, as in
section 2. We put$$ M=\sum_{i,j=1}^{n}\frac{1}{n}E_{ij}\otimes
E_{ji} +I,$$

we have$$\tilde{\omega}(M)=\sum_{i=1}^{n}E_{ii}+J=E+J,$$

hence $\tilde{\omega}(M)$ is an identity for ${\mathbb{A}}/{J}$.
Also
$$E_{lk}\cdot M=\sum_{i,j=1}^{n}E_{lk}\frac{1}{n}E_{ij}\otimes E_{ji}
+I=\sum_{i=1}^{n}\frac{1}{n}E_{li}\otimes E_{ik}
+I=\sum_{i,j=1}^{n}\frac{1}{n}E_{ij}\otimes E_{ji}E_{lk} +I=
M\cdot E_{lk},$$ for each $1\leq l,k\leq n$. Hence for each $A\in
\mathbb{A}$, we have $A\cdot M=M\cdot A$. It follow that $M$ is a
module diagonal for $\mathbb{A}$, therefore $\mathbb{A}$ is module
super-amenable by Theorem 2.6. Observe that in this case,
$J=\{0\}$, but yet $\mathbb{A}$ is not necessarily
super-amenable. This shows that the assumption that the action is
trivial from one side could not be dropped from Lemma 2.7. As a
concrete example, consider $\mathfrak G=\ell ^{1}(S)$, where
$S=[0,1]$ is a unital commutative semigroup with multiplication
$st=min \{s+t,1\}$, for $s,t\in S$, then $\mathfrak
G=\ell^{1}(S)$ and $\mathbb{A}=M_n(\mathfrak G)$ are not even
weakly amenable \cite{gro}, but still $\mathbb{A}$ is $\mathfrak
G$-module super-amenable with $J=\{0\}$.

The last example shows that there is an inverse semigroup $S$ for
which ${\ell ^{1}}(S)$ is module super-amenable but not
super-amenable. Let $(\mathbb{N}, \vee)$ be the commutative
semigroup of positive integers with maximum operation $m\vee
n=max(m,n)$, then  each element of $\mathbb{N}$ is an idempotent,
that is $E_{\mathbb N}=\mathbb N$. Hence ${\mathbb{N}}/{\approx}$
is the trivial group with one element. Therefore ${\ell ^{1}}(
\mathbb{N})$ is module super-amenable, as an ${\ell ^{1}}(
{\mathbb{N}})$-module. If ${\ell ^{1}}( \mathbb{N})$ has a
diagonal $M=\sum_{n=1}^{\infty}f_n \otimes g_n$, it should be
$M=\delta_1 \otimes \delta_1$. In this case, we have
$\delta_p\cdot M=M\cdot\delta_p\hspace{0.2cm}(p\in \mathbb{N})$,
but this equality holds if and only if, $\delta_p\otimes
\delta_1=\delta_1\otimes \delta_p$, for each $p\in \mathbb{N}$,
which is absurd. Therefore ${\ell ^{1}}( \mathbb{N})$ is not
super-amenable by \cite[Exercise 4.1.3]{run}. Note that however,
in this case, ${\ell ^{1}}( \mathbb{N})$ has an identity.

\end{document}